\documentclass[12pt,reqno]{amsart}
\usepackage[a4paper, margin=1in]
{geometry}
\usepackage[svgnames]{xcolor}
\usepackage[ruled,vlined]{algorithm2e}
\usepackage{mathtools,multicol,etoolbox}
\usepackage{makecell}
\usepackage{nameref}
\usepackage{caption}
\makeatletter
\patchcmd{\algocf@makecaption@ruled}{\hsize}{\textwidth}{}{} 
\patchcmd{\@algocf@start}{-1.5em}{0em}{}{} 
\makeatother

\SetAlCapHSkip{0pt} 
\makeatletter
\@namedef{subjclassname@2020}{%
  \textup{2020} Mathematics Subject Classification}
\makeatother

\definecolor{lightyellow}{RGB}{255, 255, 197}

\title{Classification of full del Pezzo surfaces}

\author[Streeter]{Sam Streeter}
\address{\vspace{-0.5em}School of Mathematics, University of Bristol, Woodland Road, BS8 1UG, UK}
\email{sam.streeter@bristol.ac.uk}

\usepackage{amsmath}
\usepackage{amssymb}
\usepackage{amsthm}
\usepackage{enumerate}
\usepackage{tikz-cd}
\usepackage{soul}
\usepackage[backref=page]{hyperref}
\usepackage{mathrsfs}
\usepackage{parskip}
\usepackage{booktabs}
\usepackage{tikz}
\usepackage{longtable}
\usepackage{multirow}

\theoremstyle{plain}
\newtheorem{thm}{Theorem}[section]

\newtheorem{prop}[thm]{Proposition}

 \theoremstyle{definition}
\newtheorem{defn}[thm]{Definition}

\newtheorem{rem}[thm]{Remark}

\newcommand{\thistheoremname}{}
\newtheorem{genericthm}[thm]{\thistheoremname}

\newtheorem*{genericthm*}{\thistheoremname}
\newenvironment{namedthm*}[1]
  {\renewcommand{\thistheoremname}{#1}%
   \begin{genericthm*}}
  {\end{genericthm*}}

\let\originalleft\left
\let\originalright\right
\renewcommand{\left}{\mathopen{}\mathclose\bgroup\originalleft}
\renewcommand{\right}{\aftergroup\egroup\originalright}

\DeclareMathOperator{\ch}{char}

\DeclareMathOperator{\PGL}{PGL}
\DeclareMathOperator{\Pic}{Pic}

\renewcommand{\epsilon}{\varepsilon}

\AtBeginDocument{\addtocontents{toc}{\protect\begin{multicols}{2}}}
\AtEndDocument{\addtocontents{toc}{\protect\end{multicols}}}

\begin{document}

\begin{abstract}
A split del Pezzo surface is called \emph{full} if all of its rational points lie on its exceptional curves. We complete the classification of full del Pezzo surfaces by showing non-existence in degree one and settling the remaining cases in degree two, following work of Knecht and Reyes. In degree two, $\mathbb F_{16}$ is the only field of even order admitting a full surface; we give explicit models for the three isomorphism classes over $\mathbb F_{16}$ and the two over $\mathbb F_{23}$, and we show that no such surfaces exist over larger fields of odd order.
\end{abstract}

\maketitle

\vspace{-1em}
\section{Introduction} 

This paper concerns lines on del Pezzo surfaces over finite fields. A \emph{del Pezzo surface} is a smooth projective surface $X$ over a field with ample anticanonical class $-K_X$; its degree is $d:=K_X^2 \in \{1,\dots,9\}$. Examples include smooth cubic surfaces in $\mathbb{P}^3$ (degree 3) and smooth complete intersections of two quadrics in $\mathbb{P}^4$ (degree 4). A curve $L \subset X$ is a \emph{line} (or \emph{exceptional curve}) if $L \cdot K_X = L^2 = -1$. For $d \geq 3$, the class $-K_X$ is very ample and induces an embedding in $\mathbb{P}^d$: lines are identified with lines in $\mathbb{P}^d$, hence the name.

The study of lines on del Pezzo surfaces goes back at least to the \emph{Cayley--Salmon theorem} in 1849 \cite{Cay49}, which states that a smooth cubic surface over an algebraically closed field contains exactly 27 lines. A century later, Segre \cite{Seg49} gave the full list of possible $k$-rational line counts on a cubic surface over a field $k$ and showed that they are all realised when $k = \mathbb{Q}$. Starting with Manin \cite[IV.23.13(i)]{Man86}, focus turned towards the \emph{inverse Galois problem for del Pezzo surfaces}, which asks to classify the possible Galois actions on the set of lines of del Pezzo surfaces over a given field $k$ (see e.g.\ \cite{Tre20,LT20,Karras,Trip} for finite fields and \cite{EJ15,EJ19a,EJ19b} for infinite fields). Recently, a generalisation of Segre's original result was obtained over arbitrary fields \cite{KayaMcKeanStreeterUppal}; see \cite[\S1]{KayaMcKeanStreeterUppal} for more on the history of this subject.

Over an algebraically closed field, a del Pezzo surface is isomorphic to either $\mathbb{P}^1 \times \mathbb{P}^1$ (for which $d = 8$) or the blowup of $\mathbb{P}^2$ in $9-d$ points in \emph{general position}: no 3 on a line, no 6 on a conic, and if $d = 1$, not all 8 on a cubic curve with a singularity at one \cite[IV.24.4]{Man86}. The line count is then constant for $d \neq 8$ (with two possibilities if $d = 8$) and given in the following table (8q denotes $\mathbb{P}^1 \times \mathbb{P}^1$ and 8b denotes blowups of $\mathbb{P}^2$ in a point):

\begin{longtable}{ |c|c|c|c|c|c|c|c|c|c|c| }
\caption{Geometric line counts on del Pezzo surfaces} \\
\hline
Degree & 9 & 8q & 8b & 7 & 6 & 5 & 4 & 3 & 2 & 1 \\
\hline
Number of lines & 0 & 0 & 1 & 3 & 6 & 10 & 16 & 27 & 56 & 240 \\ 
\hline
\end{longtable}

Over a general field $k$, the situation is more complicated. Del Pezzo surfaces which are isomorphic to either $\mathbb{P}^1 \times \mathbb{P}^1$ or a blowup of $\mathbb{P}^2$ in a set of $k$-rational points are known as \emph{split}. For $k$ infinite, note that split del Pezzo surfaces exist in each degree, since failure of general position corresponds to solubility of non-zero polynomials in the coordinates of the rational points (see e.g.\ \cite[Lem.~2.5]{BFL19}).

Over $k = \mathbb{F}_q$ the finite field of $q$ elements, the existence of a split degree-$d$ del Pezzo surface is more subtle. The set of $q$ for which there exists a split degree-$d$ del Pezzo surface over $\mathbb{F}_q$ is known for all $d$ thanks to work of Banwait, Fit\'e and Loughran \cite{BFL19}. 

Over finite fields, split del Pezzo surfaces can be \emph{full}, in the following sense.

\begin{defn}
A split del Pezzo surface $X/k$ is \emph{full} if every rational point $P \in X(k)$ is contained in a line of $X$.
\end{defn}

Note that no full del Pezzo surface over $k$ exists for $k$ infinite, since split $X$ is $k$-birational to $\mathbb{P}^2$, thus $X(k)$ is Zariski-dense.

Hirschfeld \cite{Hir82} classified full del Pezzo surfaces of degree $d \geq 3$, as below.

\vspace{-0.5em}
\begin{longtable}{|c|c|c|c|c|}\caption{Classification in degree $\geq 3$}\label{table:ClassAtLeast3}\\
\hline 
Degree & $\geq 6$ & 5 & 4 & 3 \\
\hline
Existence & No $q$ & $q \in \{2,3\}$ & $q =5$ & $q \in \{4,7,8,9,11,13,16\}$ \\
\hline
\#Isomorphism classes per $q$ & 0 & 1 & 1 & 1 \\
\hline
\end{longtable}
\vspace{-0.5em}
Later, Knecht and Reyes \cite{KR16} showed that full degree-2 del Pezzo surfaces over $\mathbb{F}_q$ with $q$ odd can exist only for $9 \leq q \leq 37$ and gave the following partial classification.

\vspace{-0.5em}
\begin{longtable}{|c|c|c|c|c|c|c|}
\caption{Partial classification in degree 2 for odd $q$} \label{table:Class2} \\
\hline
$q$ & 9 & 11 & 13 & 17 & 19 & 23 \\
\hline
\#Isomorphism classes & 1 & 1 & 2 & 6 & 5 & $\geq 2$ \footnote{The equation (21) given in \cite[Thm.~4.7]{KR16} is incorrect; the given surface is singular. However, the incidence data is consistent with the surface $X_2$ in Theorem~\ref{thm:main2}(2). The equation (20) in \emph{loc.\ cit.} is $X_1$ in Theorem~\ref{thm:main2}(2) (the first bitangent type listed should be $(1,6,12,3,2)$, not $(1,16,12,3,2)$).} \\
\hline
\end{longtable}
\vspace{-0.5em}
The following challenges arise in obtaining a complete classification:
\begin{enumerate}
    \item Working in even characteristic.
    \item Determining (lack of) existence for odd $25 \leq q \leq 37$ \cite[Question~4.8]{KR16}.
    \item Determining the number of isomorphism classes for $q = 23$ \cite[Question~4.9]{KR16}.
    \item Working with del Pezzo surfaces of degree $1$ \cite[Question~4.10]{KR16}.
\end{enumerate}

\subsection{Results}

Our first theorem completely settles the existence problem in degree $2$, thus partially resolving (1) and resolving (2).

\begin{thm} \label{thm:main}
\setlength{\leftmargini}{1.8em}
\begin{enumerate}
    \item There are full degree-2 del Pezzo surfaces over $\mathbb{F}_{2^n}$ if and only if $n = 4$.
    \item There are no full degree-2 del Pezzo surfaces over $\mathbb{F}_q$ for any odd $25 \leq q \leq 37$.
\end{enumerate}
\end{thm}

Our second theorem completes the resolution of (1) and resolves (3), settling the classification in degree 2.

\begin{thm} \label{thm:main2}
\setlength{\leftmargini}{1.8em}
\begin{enumerate}
\item Over $\mathbb{F}_{16}=\mathbb{F}_2(\alpha)$, where
$\alpha^4+\alpha+1=0$, there are exactly 3 isomorphism classes of full degree-2 del Pezzo surfaces, with representatives
\[
X_i:\quad w^2+A_i(x,y,z)w+B_i(x,y,z)=0
\quad\subset\mathbb{P}(1,1,1,2),
\]
where
\[
\begin{aligned}
A_1={}&x^2+xy+\alpha^4xz+\alpha^{12}y^2+\alpha^5yz+z^2,\\
B_1={}&\alpha x^3y+\alpha^5x^2y^2+\alpha^{13}xy^3
       +\alpha^6xyz^2
     +\alpha^6y^3z+\alpha^{11}y^2z^2+\alpha yz^3;
\end{aligned}
\]
\[
\begin{aligned}
A_2={}&x^2+xy+\alpha^{11}xz+y^2+\alpha^{12}yz+\alpha^4z^2,\\
B_2={}&\alpha^2x^3y+\alpha^{10}x^2y^2+\alpha^2xy^3
       +xy^2z+xyz^2
     +\alpha^6y^3z+\alpha^2y^2z^2+\alpha^5yz^3;
\end{aligned}
\]
and
\[
\begin{aligned}
A_3={}&x^2+xy+\alpha^{10}xz+\alpha^{12}y^2
       +\alpha^7yz+\alpha^{14}z^2,\\
B_3={}&\alpha x^3y+\alpha^5x^2y^2+\alpha^{13}xy^3
       +\alpha^8xy^2z+\alpha xyz^2
     +\alpha^{11}y^3z+\alpha^{14}y^2z^2+\alpha^9yz^3.
\end{aligned}
\]
\item There are exactly 2 isomorphism classes of full degree-2 del Pezzo surfaces over $\mathbb{F}_{23}$ with representatives
\[
\begin{aligned}
X_1: w^2 & = x^4 + y^4 + z^4 +6(x^2y^2 + x^2z^2 + y^2z^2), \\
X_2: w^2 & = x^4 + y^4 + z^4 +4(x^2y^2 + x^2z^2 + y^2z^2).
\end{aligned}
\]
\end{enumerate}
\end{thm}

Lastly, we make the following observation, which resolves (4) and completes the classification in all degrees.

\begin{prop}
Let $X$ be a degree-1 del Pezzo surface over any field $k$. Then the base-point $\mathcal{O} \in X(k)$ of the anticanonical linear system lies on no line of $X$. In particular, no del Pezzo surface of degree $1$ is full.
\end{prop}

\begin{proof}
Suppose $\mathcal{O}$ lies on a line $L$. Take $P \in L \setminus \{\mathcal{O}\}$. Since $h^0(X,-K_X) = 2$, there exists a unique anticanonical curve $D$ containing $P$, and $D$ is irreducible since $D \cdot (-K_X) = 1$ and $-K_X$ is ample. Then $L$ and $D$ are distinct irreducible curves with $\{P,\mathcal{O}\} \subset L \cap D$, implying $L \cdot D \geq 2$, but $L \cdot D = 1$, contradiction.
\end{proof}

\subsection{Method}

We work directly with blowups of $\mathbb{P}^2$ in seven rational points. General position and fullness for the blowup are determined by rational points on certain rational curves. Thus, finding full del Pezzo surfaces amounts to an explicit iterative procedure over 7-tuples $(P_1,\dots,P_7)$ of rational points in the plane. As we work up to isomorphism, and general position and fullness are preserved under $\text{PGL}_3$, we take $(P_1,\dots,P_4) = ([1:0:0],[0:1:0],[0:0:1],[1:1:1])$. For general position, we then have $P_i = [x_i:y_i:1]$ for $i \in \{5,6,7\}$, where $x_i,y_i \in \mathbb{F}_q \setminus \{0,1\}$ and $x_i \neq y_i$, reducing our search space to size $< q^6$. Further considering orbits under $H \leq \text{PGL}_3$ permuting the first 4 points, we reduce our fifth point candidate set to a manageable size. Explicitly,
\[
H:=\left\{\begin{pmatrix} \mathbf{v}_1^\intercal \, \vline \, \mathbf{v}_2^\intercal \, \vline \, \mathbf{v}_3^\intercal \end{pmatrix}: \mathbf{v}_i \in \{\mathbf{e}_1,\mathbf{e}_2,\mathbf{e}_3,(-1,-1,-1)\}, \, \mathbf{v}_i \neq \mathbf{v}_j\right\} \leq \text{PGL}_3,
\]
where $\mathbf{e}_1,\mathbf{e}_2,\mathbf{e}_3$ are the standard basis vectors.

We present three main algorithms:

{
\setlength{\leftmargini}{7em}
\begin{enumerate}
    \item[{\bfseries Algorithm~\ref{alg:1}:}] Produces a set of 7-tuples in $\mathbb{P}^2(\mathbb{F}_q)$ such that any full degree-2 del Pezzo surface over $\mathbb{F}_q$ is isomorphic to the blowup of $\mathbb{P}^2_{\mathbb{F}_q}$ in one of the members.
    \item[{\bfseries Algorithm~\ref{alg:2}:}] Identifies isomorphism classes of blowups of the output of Algorithm~\ref{alg:1}.
    \item[{\bfseries Algorithm~\ref{alg:3}:}] Computes anticanonical models for the output of Algorithm~\ref{alg:2}.
\end{enumerate}
}

\subsection*{Outline} \label{sec:Outline}
In Section~\ref{sec:Alg1}, we present Algorithm~\ref{alg:1}. In Section~\ref{sec:Alg2}, we present Algorithm~\ref{alg:2}, refining the output of Algorithm~\ref{alg:1} to a complete classification. In Section~\ref{sec:Alg3}, we present Algorithm~\ref{alg:3}, which computes an anticanonical model from each isomorphism class.

\subsection*{Acknowledgements}
I thank Jonathan Bober, Daniel Loughran, Ross Paterson and Happy Uppal for helpful discussions. I also thank the authors of \cite{DKPST} for discussions on how best to coordinate our work. Lastly, I thank the University of Bristol and the Heilbronn Institute for Mathematical Research for their support. All computations were performed using the \texttt{MAGMA} computer algebra package \cite{Magma}; code is available at
\vspace{-0.4em}
\begin{center}
\url{https://github.com/SamStreeter95/Full-del-Pezzos} \vspace{-0.4em}   
\end{center}

\subsection*{AI Disclosure}
The author used Claude Opus and ChatGPT-6 Astra to write/improve portions of the \texttt{MAGMA} code and in proofreading, and they retain sole responsibility for the correctness of the arguments.

\section{Algorithm 1: Partial classification} \label{sec:Alg1}

We now summarise our initial algorithm to find full del Pezzo surfaces.

\setlength{\leftmargini}{1.5em}
\begin{algorithm}[H] 
\label{alg:1}
\caption{\bf Finding full del Pezzo surfaces over $\mathbb{F}_q$}

\KwIn{A prime power $q$.}
\KwOut{A list \texttt{FullTuples} of $(a,b,c,d,e,f) \in \mathbb{F}_q^6$ such that any full degree-2 del Pezzo surface over $\mathbb{F}_q$ is isomorphic to the blowup of $\mathbb{P}^2_{\mathbb{F}_q}$ in $([1:0:0],[0:1:0],[0:0:1],[1:1:1],[a:b:1], [c:d:1], [e:f:1])$ for some $(a,b,c,d,e,f) \in \texttt{FullTuples}$.}

\begin{enumerate}
\item Create a list \texttt{AllPoints} of all elements in $\{(x,y) \in \mathbb{F}_q^2: x,y \not\in \{0,1\}, x \neq y\}$.
\item Create empty sets \texttt{FullTuples} and \texttt{Excluded}.
\item Compute the set \texttt{Orbits} of orbits of $H$ on \texttt{AllPoints} and create a set \texttt{OrbitReps} of orbit representatives for $H$ acting on \texttt{AllPoints}.
\item For each $(a,b) \in \texttt{OrbitReps}$, writing $P_5:=[a:b:1]$:
    \begin{enumerate}
    \item Set \texttt{Available}$:=$\texttt{AllPoints}$\setminus$\texttt{Excluded} and create an empty set \texttt{FiveBad}.
    \item Compute $\texttt{AllPoints} \cap L_i(\mathbb{F}_q)$ for $L_i$ the line joining $P_i$ to $P_5$ for $1 \leq i \leq 4$; add to \texttt{FiveBad}.
    \item Compute $\texttt{AllPoints} \cap C(\mathbb{F}_q)$ for $C$ the conic through $P_1,\dots,P_5$; add to \texttt{FiveBad}.
    \item Create lists $\texttt{V}:=\texttt{Available} \setminus \texttt{FiveBad}$ and $\texttt{FullV}:=\texttt{AllPoints} \setminus \texttt{FiveBad}$.
    \item For each $(c,d) \in \texttt{V}$, writing $P_6:=[c:d:1]$:
    \begin{enumerate}
        \item Create an empty list \texttt{SixBad}.
        \item Compute $\texttt{AllPoints} \cap L'_i(\mathbb{F}_q)$ for $L_i'$ the line joining $P_i$ to $P_6$ for $1 \leq i \leq 5$; add to \texttt{SixBad}.
        \item Compute $\texttt{AllPoints} \cap C_i(\mathbb{F}_q)$ for $C_i$ the conic through $P_6$ and all points in $\{P_1,\dots,P_5\}$ except $P_i$ for $1 \leq i \leq 5$; add to \texttt{SixBad}.
        \item Create lists $\texttt{W}:=\texttt{V} \setminus \texttt{SixBad}$ and $\texttt{FullW}:=\texttt{FullV} \setminus \texttt{SixBad}$.
        \item For each $(e,f) \in \texttt{W}$ after $(c,d)$ in \texttt{AllPoints}, writing $P_7:= [e:f:1]$:
        \begin{enumerate}
            \item Create an empty list \texttt{SevenBad}.
            \item Compute $\texttt{AllPoints} \cap L''_i(\mathbb{F}_q)$ for $L_i''$ the line joining $P_i$ to $P_7$ for $1 \leq i \leq 6$; add to \texttt{SevenBad}.
            \item Compute $\texttt{AllPoints} \cap C_{i,j}(\mathbb{F}_q)$ for $C_{i,j}$ the conic through $P_7$ and all the $P_k$, $k \leq 6$ except $P_i$ and $P_j$ for $1 \leq i < j \leq 6$; add to \texttt{SevenBad}.
            \item Compute $\texttt{AllPoints} \cap D_i (\mathbb{F}_q)$ for $D_i$ the cubic curve through all the $P_j$ with a singularity at $P_i$ for $1 \leq i \leq 7$; add to \texttt{SevenBad}.
            \item Define $\texttt{X}:=\texttt{FullW} \setminus \texttt{SevenBad}$.
            \item If $\#\texttt{X} = 0$, then add $(a,b,c,d,e,f)$ to \texttt{FullTuples}; else continue.
        \end{enumerate}
    \end{enumerate}
    \item Add the $H$-orbit of $(a,b)$ to \texttt{Excluded}.
    \end{enumerate}
\end{enumerate}
\end{algorithm}
\vspace{1em}
\begin{rem}
We make some observations on Algorithm~\ref{alg:1}:
\begin{enumerate}
    \item Every seven-point set in general position in $\mathbb P^2(\mathbb{F}_q)$ whose blowup is full is $\mathrm{PGL}_3(\mathbb F_q)$-equivalent to a configuration encoded by an element of \texttt{FullTuples}.
    \item The algorithm admits an immediate adaptation to finding representatives for all \emph{split} degree-2 del Pezzo surfaces over any finite field and examples over an infinite field.
    \item The calculation of line, conic and cubic equations is linear algebra, as outlined in \cite[Lem.~2.5]{BFL19}, which can be sped up by pre-computation of certain matrices.
\end{enumerate}
\end{rem}

\subsection{Proof of Theorem~\ref{thm:main}}

Running Algorithm~\ref{alg:1} for $q \in \{25,27,29,31,37\}$ we obtain $\#\texttt{FullTuples} = 0$ in each case, thus there are no full degree-2 del Pezzo surfaces over $\mathbb{F}_q$ for these $q$. It remains only to consider the case $q = 2^n$. First, note that a split degree-2 del Pezzo surface exists over $\mathbb{F}_q$ if and only if $q \geq 9$ \cite[Prop.~4.4]{BFL19}. Next, note that, for $L$ a line on a full degree-2 del Pezzo surface over $\mathbb{F}_q$, we have $\#L(\mathbb{F}_q) = q+1$ and $\#X(\mathbb{F}_q) = q^2 + 8q + 1$ (since $\#\mathbb{P}^2(\mathbb{F}_q) = q^2 + q + 1$ and each blowup point contributes $q$ points). Then for $X$ to be full, we must have $q^2 + 8q + 1 \leq 56(q+1)$, thus $q \leq 49$. It remains to consider $q \in \{16,32\}$. We find full 7-tuples over $\mathbb{F}_{16}$ and none over $\mathbb{F}_{32}$.

\section{Algorithm 2: Completing the classification} \label{sec:Alg2}

We next use \emph{Cremona transformations} to determine which blowups of the configurations returned by Algorithm 1 are isomorphic.

\begin{defn}
For $P,Q,R \in \mathbb{P}^2$ non-collinear points, the \emph{(quadratic) Cremona transformation centred at $P,Q,R$} is the birational map $\phi_{P,Q,R}:\mathbb{P}^2 \setminus \{P,Q,R\} \rightarrow \mathbb{P}^2$ defined by acting by the inverse of a matrix with columns coordinates for $P,Q,R$ (sending them to $([1:0:0],[0:1:0],[0:0:1]))$ and then applying $[x:y:z] \mapsto [yz:xz:xy]$.
\end{defn}

Given a set of points $T = \{P_1,\dots,P_n\} \subset \mathbb{P}^2$ in general position with $3 \leq n \leq 8$ and an ordered triple $(P_i,P_j,P_k)$, we obtain a new set $T' = \{Q_i,Q_j,Q_k\} \cup \{R_\ell\}_{\ell \not\in\{i,j,k\}}$ in general position, where $R_\ell$ is the image of $P_\ell$ under $\phi_{P_i,P_j,P_k}$ and $Q_i,Q_j,Q_k$ are the images of the lines $L_{P_jP_k},L_{P_iP_k},L_{P_iP_j}$ respectively. It is easily seen that the blowups at $T$ and $T'$ are isomorphic.

\begin{thm} \label{thm:iso}
Let $S$ and $S'$ be split degree-$(9-n)$ del Pezzo surfaces over a field $k$ obtained as the blowups of two size-$n$ subsets $T,T' \subset \mathbb{P}^2(k)$ with $n \leq 8$. Then $S$ and $S'$ are $k$-isomorphic if and only if $T$ and $T'$ are connected by a sequence of projective linear automorphisms and quadratic transformations centred at triples from the intermediate sets.
\end{thm}

\begin{proof}
Sufficiency is clear, so assume that $f:S \xrightarrow{\sim} S'$ is an isomorphism. The case $n \leq 4$ is clear by the transitivity of $\PGL_3$ acting on general-position $r$-tuples for $r \leq 4$, so assume henceforth that $n \geq 5$. Denote by $\pi: S \rightarrow \mathbb{P}^2$ and $\pi':S' \rightarrow \mathbb{P}^2$ the blowups in $T$ and $T'$ respectively, and let $\rho = \pi' \circ f: S \rightarrow \mathbb{P}^2$ and $\phi = \pi'\circ f \circ \pi^{-1}:\mathbb{P}^2 \dashrightarrow \mathbb{P}^2$ be the induced birational map, so that we have the following commutative diagram:
\[
\begin{tikzcd}
S \arrow[d,"\pi"] \arrow[r,"f"] \arrow[rd,"\rho"] & S' \arrow[d,"\pi'"] \\
\mathbb{P}^2 \arrow[r,dashed,"\phi"] & \mathbb{P}^2.
\end{tikzcd}
\]
Denote by $E_1,\dots,E_n$ the exceptional curves contracted by $\pi$ and set $L:=\pi^*\mathcal{O}(1)$. Then $\{L,E_1,\dots,E_n\}$ is a basis of $\Pic(\overline{S})$. Set $D:=\rho^*\mathcal{O}(1)$. Then we may write 
\[
D = dL - \sum_{i=1}^nm_i E_i,
\]
where $d$ is the degree of $\phi$ and each $m_i \geq 0$ since $D$ is nef and $D \cdot E_i = m_i$. Note that
\[
D^2 = d^2 - \sum_{i=1}^nm_i^2 = 1, \quad -K_S \cdot D = 3d - \sum_{i=1}^n m_i = 3,
\]
thus we obtain
\[
\sum_{i=1}^nm_i^2 = d^2 -1, \quad \sum_{i=1}^nm_i = 3d-3.
\]
Without loss of generality, we may assume that $m_1 \geq m_2 \geq \cdots \geq m_n$. Then 
\begin{equation} \label{eq:ineq}
\begin{aligned}
d^2 -1 & = \sum_{i=1}^nm_i^2 = m_3\sum_{i=1}^nm_i + \sum_{i=1}^n m_i(m_i - m_3) \\
& \leq 3(d-1)m_3 + (d-1)(m_1 - m_3) + (d-1)(m_2 - m_3) \\
& = (d-1)(m_1 + m_2 + m_3), 
\end{aligned}
\end{equation}
where we have used in the inequality that each $m_i \leq d-1$, since $\sum_{i=1}^nm_i^2 = d^2-1$.

Suppose that $d > 1$. Then we deduce from Equation~\eqref{eq:ineq} that $m_1 + m_2 + m_3 > d$. Denote by $\psi: \mathbb{P}^2 \dashrightarrow \mathbb{P}^2$ the quadratic Cremona transformation centred at the three points $P_i:=\pi(E_i)$ for $i \in \{1,2,3\}$ and by $\pi_1: S \rightarrow \mathbb{P}^2$ the morphism extending $\psi \circ \pi$. We obtain a new commutative diagram:
\[
\begin{tikzcd}
S \arrow[d,"\pi_1"] \arrow[r,"f"] \arrow[rd,"\rho"] & S' \arrow[d,"\pi'"] \\
\mathbb{P}^2 \arrow[swap,r,dashed,"\phi \circ \psi^{-1}"] & \mathbb{P}^2.
\end{tikzcd}
\]
We have $L_1:=\pi_1^*\mathcal{O}(1) = 2L - E_1 - E_2 - E_3$, and the degree of $\phi \circ \psi^{-1}$ is $d_1:=D\cdot L_1 = 2d - (m_1 + m_2 + m_3) < d$, strictly smaller than $d$.

Repeating this process while the degree of the birational map is greater than $1$, we must eventually reduce it to degree 1, in which case it is a projective linear automorphism carrying the final set of blowup centres to $T'$.
\end{proof}

Consequently, we may determine isomorphism classes of full del Pezzo surfaces from the output of Algorithm~\ref{alg:1}.

One may assign a canonical $\PGL_3$-representative $\texttt{Normal}(T)$ to each 7-tuple $T$ from the output of Algorithm~\ref{alg:1} using the following algorithm (for which we fix a total order on $\mathbb{F}_q$).

\begingroup
\renewcommand{\thealgocf}{N}
\addtocounter{algocf}{-1}
\begin{algorithm}[H]
\label{alg:N}
\caption{\bf Creating canonical $\text{PGL}_3$-orbit representatives} 
\KwIn{A general position set $T = \{P_1,\dots,P_7\} \subset \mathbb{P}^2$.}
\KwOut{A representative $\texttt{Normal}(T)$ for the $\text{PGL}_3$-equivalence class of the input.}

\begin{enumerate}
    \item Construct a function \texttt{KeyPt} which canonically normalises any projective point by dividing all coordinates by the first non-zero coordinate.
    \item Initialise \texttt{Best} as empty and set \texttt{Standard}$:=([1:0:0],[0:1:0],[0:0:1],[1:1:1])$.
    \item For each of the ordered 4-tuples $S:=(P_i,P_j,P_k,P_l)$ of distinct points from $T$, apply to $T$ the matrix $M_S \in \text{PGL}_3$ sending $S$ to \texttt{Standard}. Normalise (via \texttt{KeyPt}) and sort the three remaining image points $(Q_5,Q_6,Q_7)$ lexicographically. If \texttt{Best} is empty or $(Q_5,Q_6,Q_7)$ is less than \texttt{Best}, set \texttt{Best}$:=(Q_5,Q_6,Q_7)$.
    \item Return \texttt{Normal(T)}$:=$concatenation of \texttt{Standard} and \texttt{Best}.
\end{enumerate}
\end{algorithm}
\endgroup
\vspace{-0.5em}
The next algorithm then computes isomorphism classes after Algorithm~\ref{alg:1}. \\
\begin{algorithm}[H]
\caption{\bf Checking for $\mathrm{PGL}_3$ and Cremona equivalence}
\label{alg:2}

\KwIn{The list $L$ of configurations produced by Algorithm~\ref{alg:1}.}
\KwOut{A list \texttt{CremonaReps} of representatives up to
$\mathrm{PGL}_3$ and Cremona equivalence.}

\begin{enumerate}
    \item Compute
    $\texttt{Reps}:=\{\texttt{Normal}(T):T\in L\}$,
    listed as $R_1,\dots,R_n$.
    \item Create an associative array \texttt{ClassIndex},
    setting $\texttt{ClassIndex}[R_i]:=i$ for each $i$.
    \item For each $i\in\{1,\dots,n\}$:
    \begin{enumerate}
        \item Initialise $\texttt{Neighbours}[i]$ as an empty set.
        \item For each of the $35$ unordered triples $S$ from $R_i$:
        \begin{enumerate}
            \item Apply the Cremona transformation centred at $S$.
            Form $T$ from the images of the four remaining points
            and the three points obtained by contracting the lines
            joining pairs of points in $S$.
            \item Set $R:=\texttt{Normal}(T)$ and
            $j:=\texttt{ClassIndex}[R]$.
            Add $j$ to $\texttt{Neighbours}[i]$.
        \end{enumerate}
    \end{enumerate}
    \item Initialise \texttt{Visited} as an empty set and
    \texttt{CremonaReps} as an empty list.
    \item For each $i\in\{1,\dots,n\}$ not already in \texttt{Visited}:
    \begin{enumerate}
        \item Append $R_i$ to \texttt{CremonaReps}
        and set $\texttt{Pending}:=\{i\}$.
        \item While \texttt{Pending} is non-empty:
        \begin{enumerate}
            \item Remove an index $j$ from \texttt{Pending}.
            \item If $j\notin\texttt{Visited}$, add $j$ to
            \texttt{Visited} and add all members of
            $\texttt{Neighbours}[j]\setminus\texttt{Visited}$
            to \texttt{Pending}.
        \end{enumerate}
    \end{enumerate}
    \item Return \texttt{CremonaReps}.
\end{enumerate}
\end{algorithm}
\vspace{-0.75em}
\section{Algorithm 3: Computing anticanonical models} \label{sec:Alg3}
Over a field $k$ of characteristic $\neq 2$, any degree-2 del Pezzo surface $X/k$ is of the form
\[
X: w^2 = A(x,y,z) \subset \mathbb{P}(1,1,1,2),
\]
where $A(x,y,z)$ is homogeneous of degree $4$.

When $\ch k = 2$, any such $X/k$ instead has an Artin--Schreier equation
\[
X:w^2 + wA(x,y,z) + B(x,y,z) = 0 \subset \mathbb{P}(1,1,1,2),
\]
where $A$ and $B$ are homogeneous of degrees 2 and 4 respectively.

The branch curve of the anticanonical morphism $X \rightarrow \mathbb{P}^2$ is, in either case, $A(x,y,z) = 0$.

To determine the branch curve from a 7-tuple, we consider the following diagram:

\begin{center}
\begin{tikzcd}
X \arrow[d,swap,"\rho"] \arrow[dr,"\phi"] & \\
\mathbb{P}^2 \arrow[r,dashed,"\psi"] & \mathbb{P}^2. 
\end{tikzcd}
\end{center}

Here, $\rho: X \rightarrow \mathbb{P}^2$ is the blowup of $\mathbb{P}^2$ in the 7-tuple, $\psi: \mathbb{P}^2 \dashrightarrow \mathbb{P}^2$ is the rational map induced by a basis $\{F_0,F_1,F_2\} \subset \mathcal{O}_{\mathbb{P}^2}(3)(\mathbb{P}^2)$ for the linear system of cubic curves through the 7-tuple, and $\phi$ is the resolution of indeterminacy of $\psi$. The Jacobian curve $J(x,y,z) = 0$ (the locus of singular points of singular members of the system) is the plane image of the anticanonical ramification divisor \cite[\S8.7.2]{Dolgachev2012}. Then it suffices, for each full 7-tuple we have found, to calculate a basis $\{F_0,F_1,F_2\}$ and use it to compute $\psi$ and $J$. Let $M(x,y,z)$ be the $3 \times 3$ matrix of partial derivatives of the $F_i$. For $\ch(k) \neq 3$, we have $J(x,y,z)=\det M$, while for $\ch(k) = 3$, we have $J(x,y,z)=1/z \det N$, where $N(x,y,z)$ is $M$ with the row of $z$-partial derivatives replaced by $(F_0,F_1,F_2)$. Write $F=(F_0,F_1,F_2)$, and choose a sextic $S(x,y,z)$ vanishing to order at least two at each point of the 7-tuple and outside the span of the products $F_iF_j$. Then $S$ satisfies some quadratic relation $S^2 + a(F)S + b(F) = 0$ with $a$ and $b$ homogeneous of degrees 2 and 4 respectively. Its $S$-derivative $2S+a(F)$ and $J$ both cut out the ramification divisor of $\psi$ away from the 7-tuple, hence $2S + a(F) = \mu J$ for some $\mu \in k^\times$. Completing the square in odd characteristic and rescaling in even characteristic, we obtain
\[
A(F_0,F_1,F_2) =
\begin{cases}
J(x,y,z)^2 \text{ if $\ch k \neq 2$}, \\
J(x,y,z) \text{ if $\ch k = 2$}.
\end{cases}
\]
Recovering $A$ is then a $k$-linear algebra problem in its coefficients.

In the case $\ch k = 2$, we must find a homogeneous sextic $S$ not in the $k$-span of the products $F_iF_j$, $0 \leq i \leq j \leq 2$, a homogeneous quartic $B_0$ and scalar $\mu \in k^\times$ such that
\[
S(x,y,z)^2 + \mu A(F_0,F_1,F_2)S(x,y,z) + B_0(F_0,F_1,F_2) = 0.
\]
Setting $B:=\mu^{-2}B_0$, the anticanonical model is
\[
w^2 + wA(x,y,z) + B(x,y,z) = 0.
\]

To do this, we first compute a basis for the 7-dimensional space of sextics with a double point at each of the blowup points $P_i$ and choose a generator $S(x,y,z)$ for the complement of the 6-dimensional subspace generated by the products $F_iF_j$. We then solve linearly for the 15 unknown coefficients of $B_0(x,y,z)$ and the scalar $\mu$.

\begin{algorithm}[H] 
\caption{\bf Computing anticanonical models over a field $k$}
\label{alg:3}

\KwIn{A 6-tuple $(a,b,c,d,e,f) \in k^6$ such that the following set is in general position: $\{[1:0:0],[0:1:0],[0:0:1],[1:1:1],[a:b:1],[c:d:1],[e:f:1]\} \subset \mathbb{P}^2(k)$.}
\KwOut{An anticanonical model for the blowup of this 7-tuple.}

\begin{enumerate}
\item Compute a basis $\{F_0,F_1,F_2\}$ for the system of cubics through the input set.
\item Compute the Jacobian sextic $J$.
\item If $\ch k \neq 2$, then solve for the 15 coefficients of $A(x,y,z)$ such that $J^2 = A(F_0,F_1,F_2)$ and return $w^2 = A(x,y,z)$. Otherwise, solve for the 6 coefficients of $A$ such that $J = A(F_0,F_1,F_2)$, then the $28$ coefficients of a sextic $S(x,y,z)$ with a double point at each element of the 7-tuple but not in the span of the products $F_iF_j$, then for the $15$ coefficients of $B_0(x,y,z)$ and a scalar $\mu$ such that $S^2 + \mu JS + B_0(F_0,F_1,F_2) = 0$ and return $w^2 + wA + B = 0$ with $B = \mu^{-2}B_0$.
\end{enumerate}
\end{algorithm}

\begin{proof}[Proof~of~Theorem~\ref{thm:main2}]
Running Algorithm~\ref{alg:2}, we obtain the listed numbers of classes. Running Algorithm~\ref{alg:3}, we obtain models isomorphic to the listed ones (verified in accompanying \texttt{MAGMA} code).    
\end{proof}

\bibliography{del-pezzo}{}
\bibliographystyle{alpha}

\end{document}